\documentclass[journal]{IEEEtran}

\usepackage{fixltx2e}
\usepackage{url}
\usepackage{multirow}
\usepackage{float}
\usepackage{graphicx}
\usepackage{amssymb,amsmath}
\usepackage{textcomp}
\usepackage{mathtools}
\usepackage{xfrac}
\usepackage{subfig}
\usepackage{epstopdf}

\usepackage{xfrac}
\usepackage{multirow}
\usepackage{fixltx2e}
\usepackage{url}
\usepackage{multirow}
\usepackage{wrapfig}
\usepackage{booktabs}
\usepackage{amssymb}
\usepackage{rotating}
\usepackage{array}
\usepackage{bigstrut}
\usepackage{paralist, amssymb}
\usepackage{lipsum}

\newcommand{\bc}{,\;\penalty0 }

\makeatletter
\newcommand*{\rom}[1]{\expandafter\@slowromancap\romannumeral #1@}
\makeatother
\ifCLASSINFOpdf
\else
\fi
\begin{document}
%
\title{Electric Power Allocation in a Network of Fast Charging Stations}
\author{I.~Safak~Bayram,~\IEEEmembership{Graduate~Student~Member,~IEEE,} George~Michailidis,~\IEEEmembership{Member,~IEEE, }
Michael~Devetsikiotis,~\IEEEmembership{Fellow,~IEEE,} and~Fabrizio~Granelli~\IEEEmembership{Senior~Member,~IEEE}
\thanks{Manuscript received 8 October 2012; revised 18 March 2013. Part of the work appeared in~\cite{sgc12} at Smartgridcomm'12 Conference, Tainan City, Taiwan.}
\thanks{I.~Safak Bayram and Michael Devetsikiotis are with the Department
of Electrical and Computer Engineering, North Carolina State University, Raleigh,
NC, 27695-7911 USA.}
\thanks{George Michailidis (corresponding author)
is with the Departments of Statistics and EECS, University of Michigan, Ann Arbor, MI, 48109-1107, USA.}
 \thanks{Fabrizio Granelli is with the Department of Information Engineering and Computer Science, University of Trento, Trento, Italy.}
  \thanks{Emails:(isbayram, mdevets)@ncsu.edu, gmichail@umich.edu and granelli@disi.unitn.it.}}

\maketitle
\begin{abstract}
In order to increase the penetration of electric vehicles, a network of fast
charging stations that can provide drivers with a certain level of quality
of service (QoS) is needed. However, given the strain that such a network can
exert on the power grid, and the mobility of loads represented by electric
vehicles, operating it efficiently is a challenging and complex problem.
In this paper, we examine a network of charging stations equipped with an
energy storage device and propose a scheme that allocates power to them
from the grid, as well as routes customers. We examine three scenarios,
gradually increasing their complexity. In the first one, all stations
have identical charging capabilities and energy storage devices,
draw constant power from the grid and no routing decisions of customers
are considered. It represents the current state of affairs and serves as
a baseline for evaluating the performance of the proposed scheme. In the
second scenario, power to the stations is allocated in an optimal manner
from the grid and in addition a certain percentage of customers can be
routed to nearby stations. In the final scenario, optimal allocation of
both power from the grid and customers to stations is considered.
The three scenarios are evaluated using real traffic traces corresponding
to weekday rush hour from a large metropolitan area in the US. The results
indicate that the proposed scheme offers substantial improvements of
performance compared to the current mode of operation; namely,
more customers can be served with the same amount of power, thus enabling
the station operators to increase their profitability. Further, the scheme
provides guarantees to customers in terms of the probability of being
blocked (and hence not served) by the closest charging station to their location.
Overall, the paper addresses key issues related to the efficient operation, both from
the perspective of the power grid and the drivers satisfaction, of a network
of charging stations.
\end{abstract}

\noindent\begin{keywords}
Electric Vehicles, Stochastic Charging Station Model, Performance Evaluation
\end{keywords}

\section{Introduction}\label{intro}

Over the last few years a strong push is occurring to reduce the use of hydrocarbons in transportation.
This trend is supported by the latest advances in battery and converter technology, along with
government mandates on energy independence and resilience and is enabled by the introduction of
electric vehicles (EVs) and their close relatives Plug-in Hybrid Electric Vehicles (PHEVs) by major car
manufacturers that have drastically increased consumer choices \cite{2,3}. Although there are diverging
forecasts about the growth rate of the EV population \cite{ieaRep}, there is consensus that it is going to represent
a sizable portion of the US fleet by 2025 - 30. Obviously, penetration
rates could be significantly higher than these estimates depending on battery costs, gasoline prices,
government policies, and the availability of {\em charging infrastructure}.

Indeed, such infrastructure is mostly needed in metropolitan areas, primarily characterized by higher population density,
and where residents living in multi-unit dwellings do not have easy access to night-time charging capabilities.
A recent survey among EV drivers in California shows that 40\% of them travel daily farther than the
range of their fully charged battery \cite{transDavis}, thus requiring a recharge during daytime operation of the vehicle. A network of fast charging
stations overcomes this problem \cite{ieaRep}.

On the other hand, there is concern about the strain that a rapid adoption of EVs would exert on the power grid,
due to the large load that they represent \cite{mit}. Obviously, the extent of their impact will depend on the
degree and local/regional density of the EV penetration rate, charging requirement and the time of the day
they are charged. Nevertheless, deploying large scale charging stations may lead to
grid instabilities. However, equipping each station with an energy storage device can reduce the impact of EV charging
as shown in \cite{bayramSpringer,sgc13}.

The previous discussion indicates that efficient operational regimes for a network of charging stations need to be developed,
so that they minimize the strain on the power grid, while at the same time offering good quality of service to EV drivers.
The aim of this study is to address these issues in a comprehensive manner. Specifically,
\begin{itemize}
   \item We introduce an EV fast DC charging station architecture, introduce a stochastic model to capture its operational
   characteristics and evaluate its performance (defined as the percentage of served customers). The charging station is equipped with
   a local energy storage device that aids smoothing the stochastic customer demand.
   \item We propose a resource allocation framework that meets \emph{QoS} targets at each station and
   minimizes the amount of power employed. This framework is evaluated under three different scenarios motivated by examining
   actual traffic traces from the Seattle area that exhibit a non-uniform spatial distribution of vehicles trips. The three scenarios
   in increased complexity are:
        \begin{inparaenum}[(i)]
\item no power or customer allocation to station occurs, the stations are identical in nature and act as inert service points;
\item power resources are allocated to each station; and
\item optimal power resource allocation is complemented by customer rerouting to neighboring stations.
\end{inparaenum}
   \item To achieve allocation of customers, a two-way communications protocol is introduced that coordinates EV
   assignments and reroutes if necessary. The latter represents an essential element for increasing the number of EVs being charged with
   the same amount of power drawn from the grid.
 \end{itemize}

The remainder of this paper is organized as follows: Section~\ref{prework} discusses related literature. In section~\ref{singleStation}, we introduce our single charging station stochastic model, while in section~\ref{chargingNetwork}, we introduce the network and present our power resource allocation framework for the aforementioned cases. In section~\ref{results}, we analyze the Seattle traffic traces to estimate
traffic arrival rates for different stations in the area under consideration. Further, we employ ideas from response surface methodology
to estimate a model of the performance metric of interest -probability of an EV arriving to a charging station and being blocked from
receiving service- as a function of the station's charging capacity, speed of battery charging and vehicles arrival rate.
 The resulting model is used to solve the various allocation problems introduced in section~\ref{chargingNetwork} and their results
 are compared.

\section{Related Work}\label{prework}
There has been increased interest on devising schemes that efficiently schedule EV chargings, on developing architectures for
charging station, and for organizing and operating a network of charging stations. The following paragraphs provide a brief overview of related literature.

Most works on scheduling EV chargings, assume stationary vehicles located at customer premises or large parking lots. The proposed charging strategies can be classified into the following two categories. In the first one, there is a central authority (dispatcher)
that to a large extent controls and mandates charging rates, start times, etc.
\cite{langTong, centralControl1, alizadeh2012packet, centralControl4}.
System level decisions involve selecting the desired state of charge, charging intervals, etc. are taken so as to finish all charging requests by a prespecified deadline (e.g. 7 am). The main
advantage of a centrally controlled charging schedule is that it leads to higher utilization of grid resources, together with real
time monitoring of operational conditions across the entire power system. The second category examines decentralized decision making by
EV owners. Specifically, they select individual charging patterns based on the prevailing price of electricity or on self-imposed deadlines.
It eliminates the need for a third party controller (dispatcher) and complex monitoring techniques. Since decisions are taken individually, game theoretic models, such as mean field games, potential games, and network routing games are used in these studies
~\cite{comparison2, decentralized2, decentralized3, decentralized4}.

As will be seen in section~\ref{results}, we consider spatially distributed, our study uses a centralized decision making mechanism
for a subset of EVs.

Currently there are only a handful of studies on charging station design. From a pure power engineering perspective~\cite{10} proposes a fast charging station architecture
with a DC bus distribution system. The station is equipped with an energy storage unit to minimize the strain on the grid, and the sizing problem was determined
by Monte Carlo simulations according to average load. A similar station architecture was used in ~\cite{mcgill,mcgill2}, but two different energy storage devices were considered; a flywheel and a supercapacitor. A mechanism that simultaneously draws power from the grid
and the storage devices was introduced to decrease the EVs charging duration. However, there is a multitude of storage technologies in the market
and the choice of the most appropriate one is mostly station dependent (e.g. a low energy density, large size but inexpensive storage
device may not be suitable for a station located at or near city centers, due to real-estate costs)~\cite{sgc12,globecom}. Thus, in our station architecture
we examine different storage technologies, characterized by their efficiencies and power ratings.

\section{Charging Station Architecture}\label{singleStation}

The design of a network of charging stations is ultimately linked to the current power grid operations.
At present, customer demand -household, commercial and industrial- can be assigned to three categories based on service costs.
The first represents the base load that is supplied by large, low cost (per kWh) generation assets, such as nuclear, coal and hydro.
Large size industrial customers with fairly steady demand, together with an aggregate estimate of
households and commercial users belong to this category. The second category represents the difference between base load generation and
expected aggregate demand and is primarily met by gas/liquid fuel power stations. Finally, the third category represents peak demand
that is met by fast start generators, which are characterized by their high cost (per kWh).

EVs represent sizable, {\em mobile} electric loads. Level-1 charging represents a load comparable to a household, while Level-2 charging
a load twice as large to a household. Thus, large number of EVs, geographically concentrated, would impose huge strains not only on power
generation, but also on the grid distribution system~\cite{mit,9,8}. In some studies
\cite{austin,oakridge} it is argued that if just $5$\% of all EVs charge simultaneously at fast charging stations, $5$ GW of extra power would be needed by year $2018$ in the VACAR region (Virginia - North Carolina - South Carolina). For these reasons, charging station designs that do not stress the power grid and eliminate the need for adding significant extra generation capacity become important.

\begin{figure}[t]
  \centering
  \includegraphics[width=0.8\columnwidth]{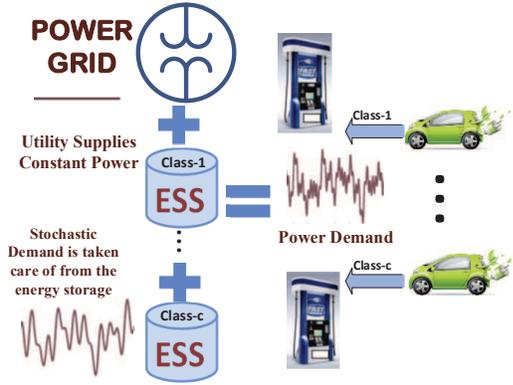}
  \caption{Single Charging Station Architecture}\label{chargingArc}
\end{figure}

The four key components of our design that try to address these issues are: (I) each station draws constant power from the grid; (II) local energy storage is employed to meet stochastic customer demand; (III) the station supports different classes of charging requests (fast service vs
slow service); and (IV) the QoS metric employed is the long-term blocking probability of incoming customers.
The overview of the proposed charging station is depicted in Figure~\ref{chargingArc}. Next, we explain the system dynamics in detail.

(I) Charging stations of any significant size represent commercial size loads. Hence, it seems reasonable for station operators to draw
long-term contracts with the utility where a power level is agreed in return for a lower price. This enables the utility to better anticipate
its demand, and the station operator to benefit from a lower price; as argued in \cite{seanMeyn, arellano2007model},such contracts leads to lower contract, as well as average spot prices, and more efficient market equilibria.

\begin{figure}[t]
 \centering
 \includegraphics[width=\columnwidth]{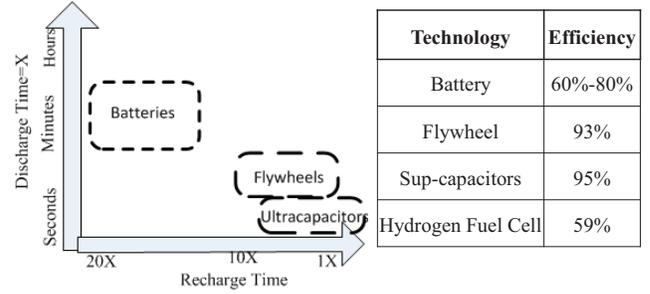}
  \caption{Candidate Energy Storage Systems Technology Landscape~\cite{essTech}~\cite{sandia}}\label{rateComp}
\end{figure}

(II) Energy storage represents a critical component in the proposed system architecture, since it aids in smoothing customers' stochastic
demand. During rush hour, stored energy can be used to serve more customers. Similarly, when the power that the station can draw from
the grid is not fully utilized, the extra power can recharge the storage device. An overview of candidate technologies for storage
devices and their efficiencies are presented in Figure~\ref{rateComp}. Their details will be further explained in the next section. (\ref{queueingModel}).

(III) We consider that the charging station provides service to multiple customer classes at different charging rates.
This allows the station to accommodate customers with different charging needs and preferences, as well as
EVs with different technological constraints.

(IV) As discussed in the introduction, charging times depend on the level used, but on average they are about 30 minutes.
Hence, it is reasonable to assume that incoming customers would not be willing to wait and thus in our model (discussed next)
a ``bufferless" system was adopted. For such a system, the blocking probability becomes a natural performance metric.

\subsection{Stochastic Model for Station Dynamics}\label{queueingModel}

Based on the aforementioned specifications, the proposed station architecture and the corresponding
model for its behavior over time, exhibit the following operation characteristics:
\begin{inparaenum}[(i)]
\item the charging station draws a constant power from the grid;
\item upon exceeding the available grid power, the local energy storage unit is used to meet additional demand;
\item whenever there is idle grid power, it is used to charge the local energy storage device, if it is not in a fully charged state;
\item depending on the amount of constantly drawn grid power and the size of the local energy storage, a certain level of \emph{QoS} is provided; and
\item the station partitions its capacity with respect to demand for each customer class. Such insights can be obtained from profiling studies (e.g. customer surveys, etc).
\end{inparaenum}
The constantly drawn grid power is discretized to $S$ equal slots, meaning that it can accommodate up to $S$ vehicles at the same time. In a similar way, the local energy storage can charge $R$ vehicles in a fully charged state. Since the charging station can never serve more than $S+R$ vehicles concurrently, the very next EV arrival is going to be ``blocked". This strategy insulates to a large extent the power grid from peak demand.
\begin{figure*}[t]
 \centering
 \includegraphics[scale=0.45]{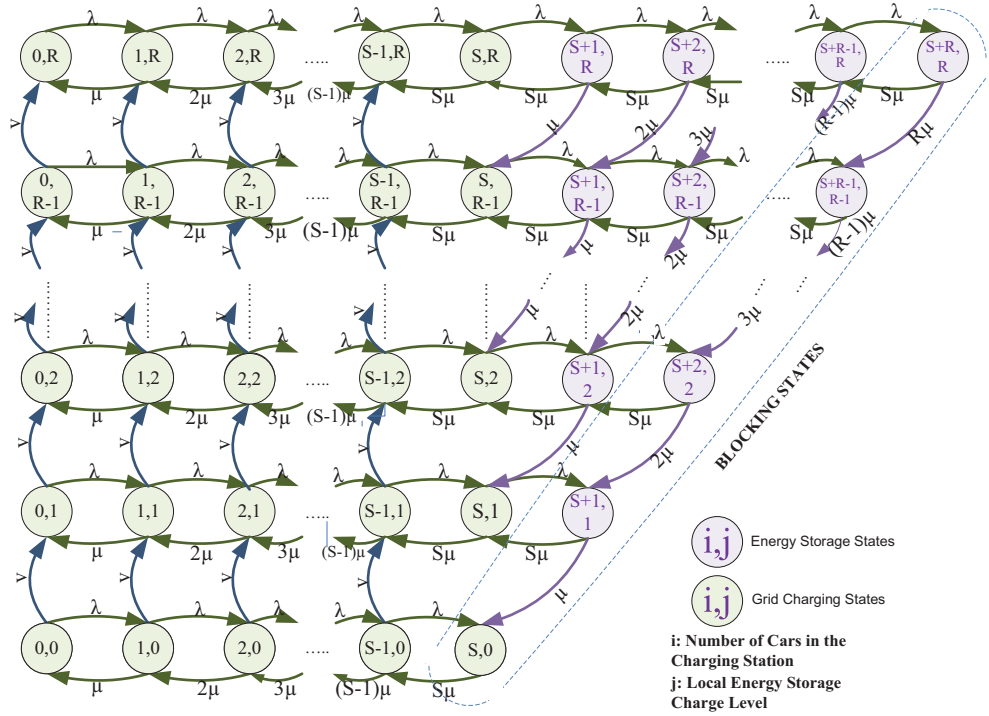}
  \caption{Continuous Time Markov Chain}\label{markovChain}
\end{figure*}
The details of the stochastic model are given next. Customers arrive to the charging facility according to a Poisson process with parameter $\lambda$. Currently, a variety of different EV models with different battery sizes exist. Thus, the service time of customers is assumed to be a exponentially distributed with rate $\mu$. Also, the charging duration of the energy storage device so that it is able to accommodate one more
EV is exponentially distributed with rate $\nu$ which depends on the underlying energy storage technology.

We proceed to explain the energy storage component in detail. The power rating of an energy storage device determines how fast it can be charged. Thus, energy stored in unit time can be calculated by the product of the power rating and the efficiency (ratio of stored energy and total amount of energy spent to charge battery) of the energy storage unit. This means that depending on these two parameters, different amounts of energy can be stored. For instance, assume that our fast charging station can charge an EV (battery with $\eta=0.9$) in $30$ minutes using maximum power rating, $S_{PR}=1$ ($\mu=2$). Also, we employ an energy storage with the same efficiency, but with a higher power rating $S_{PR}=2$ than the EV battery, so that in the same amount of time, we can store up the energy for the demand coming
from 2 EVs in the local storage device ($\nu=4$). Note that the charging rate is $\tilde\nu= f(\hat{S},\eta)$, with $\hat{S}\leq\;S_{PR}$
being the available power.

Given the assumptions above, the single charging station model can be represented by continuous time Markov chain with 2-dimensional finite state
space. In Figure~\ref{markovChain}, the state space of the Markov chain, along with its transmission rates are depicted. The total number of states is $\kappa =  \emph{(S+1)(R+1)}+\sum\limits_{i=1}^{\emph{R}} i$. It is easy to see that a unique steady state distribution would exist which
can be calculated by solving:
\begin{equation}\label{generator1}
    \pi\emph{Q} = 0  \textbf{ } \textrm{ and }\textbf{ }\pi\emph{e} = 1
   \end{equation}
where \emph{e} is a column vector whose elements are all equal to 1, and $Q$ is a $\kappa\times\kappa$ matrix containing the
transition rates and $\pi$ a vector of length $\kappa$ containing the steady state probabilities.
Note that the elements of $Q$ satisfy $q_{ab} \geq 0$ for $a \neq b $ and $q_{ab} = -\sum_{a=1,b\neq a}^\kappa q_{ab}$ for all $a = 1, 2, \ldots,\kappa$. Then, the model's (station's) blocking probability can be calculated from
$\sum\limits_{i=1}^{S} \pi(\frac{i(i+2S+1)}{2}), i=1\dots S$.

\begin{equation}\label{Qgenerator}
 Q =
 \begin{pmatrix}
 -(\lambda+\nu) & \lambda & \cdots & 0 \\
  \mu & -(\lambda+\nu+\mu) & \cdots & 0 \\
  \vdots  & \vdots  & \ddots & \vdots  \\
  0 & 0 & \cdots & -(S+R)\mu
 \end{pmatrix}
\end{equation}

Next, we extend this model to the case where different classes of customers are present; namely, $c\in \{1,\cdots, C\}$. Also, denote by
$\vec{\rho}$ the percentage of customers that demand class-$c$ type of service. Then, the station operator partitions the power drawn from the
grid into $C$ components by solving the following optimization problem.
\begin{equation}\label{multiClass}
\begin{aligned}
& arg \mathop {\min }\limits_{{S^{(c)}}}  & & \sum\limits_{c \in \mathcal{C}} {{B^{(c)}}} (\vec{\rho}\lambda , S^{(c)},  R^{(c)}) &\\
& \text{s.t. } & & \sum\limits_{c\in C} S^{(c)} = S \\
&  & &{\vec R}^{(c)},\;\vec{\rho},\; and \;\lambda\; are\; given \\
\end{aligned}
\end{equation}

\begin{figure}[t]
  \includegraphics[width=\columnwidth]{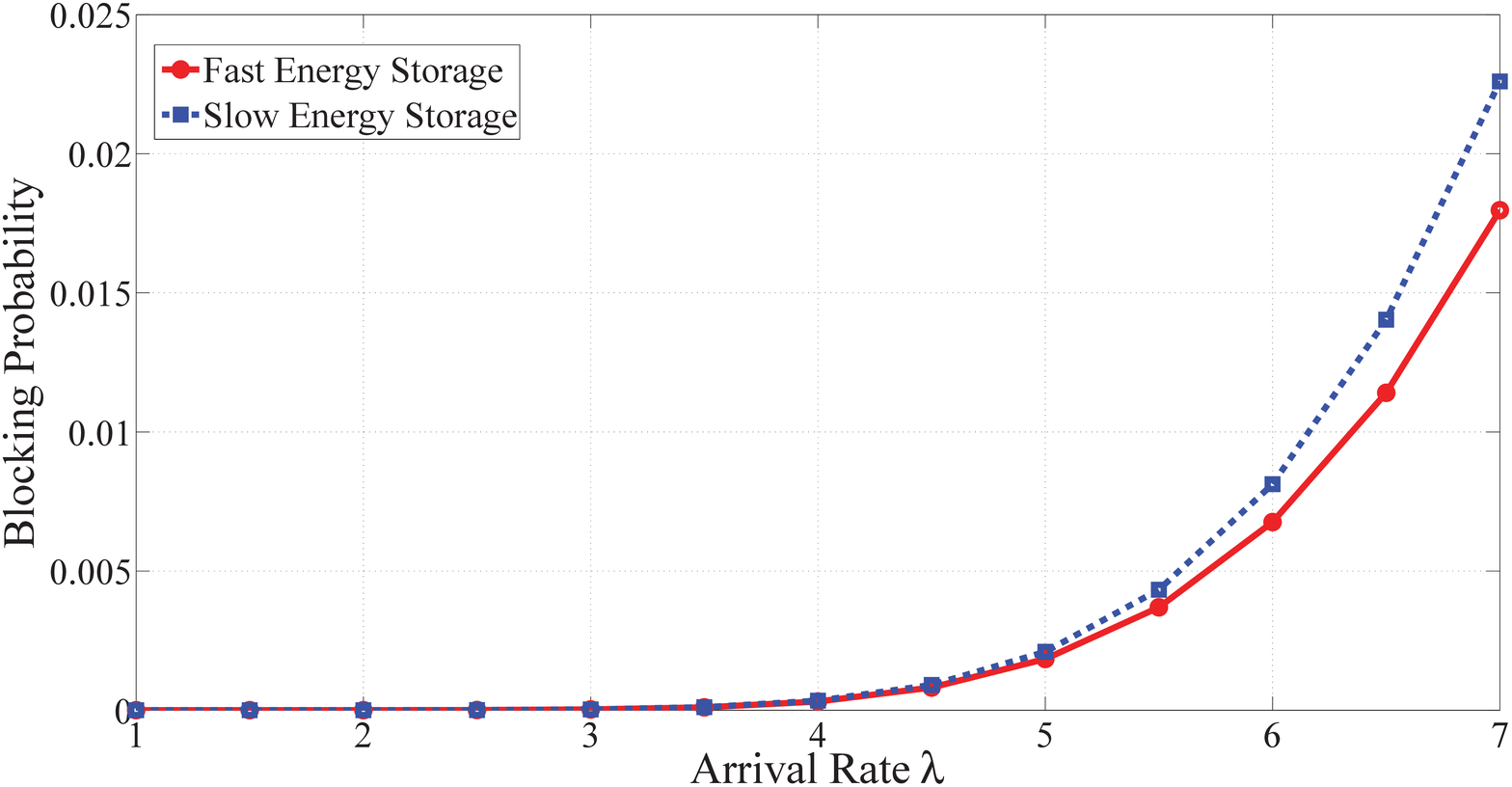}
  \caption{Performance Evaluation of Different Energy Storage Devices}\label{slowFast}
\end{figure}

To illustrate how the characteristics of the energy storage device improve performance of the station we use the following example.
We fix the size of two devices, but we vary their efficiency and power rating parameters. There is a fast energy storage with $95\%$ efficiency and $S_{PR}$=$2$, and a slow one with efficiency of $85\%$ and $S_{PR}$=$1$. Storage size is set to $R$=$5$ and the EV arrival rate varies between ($\lambda=1 - 7$). To ease the demonstration, a single customer class is assumed requesting a charging rate of
$\mu$=$2$. As shown in Figure~\ref{slowFast}, the fast energy storage device outperforms the slow one in terms of blocking probabilities.

\begin{figure}[t]
  \includegraphics[width=\columnwidth]{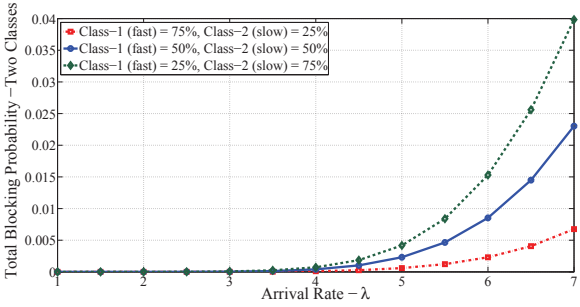}
  \caption{Multi Class Charging Station Performance Evaluation}\label{perfEva}
\end{figure}

Next, we evaluate the system performance (percentage of vehicles it can charge), under the following sets of parameters.
There are two customer classes;  in class-1 EVs request fast charging, while in class-2 request slower charging.
A typical charging duration takes $30$ minutes, thus the charging rate $\mu^{(1)}$ is set to $2$ and $\mu^{(2)}$=$1$.
We assume that the station operator picks the energy storage according to the following specifications (note that superscript denotes the customer class): storage size $R^{(1)}$=$R^{(2)}$=$5$, efficiency $\eta^{(1)}$=$0.95$ and $\eta^{(2)}$=$0.85$ and power ratings $S_{PR}^{(1)}=2$ and $S_{PR}^{(2)}$=$1$. Based on an EV profiling study, it is estimated that the total arrival rates varies between $\lambda$ = $1-7$. We look at
three different compositions of the EV population: $(\rho^1,\rho^2)=\left\{ {(75\% ,25\% ),(50\% ,50\% ),(25\% ,75\% )} \right\}$. Then,
the station operator solves optimization problem~\ref{multiClass} to calculate the optimal $\vec{S}^{(c)}$ given by  $[6,\;4]$, $[4,\;6]$, and $[2,\;8]$ for the given $(\rho^1,\rho^2)$ pairs, respectively. The resulting blocking probabilities are shown in Figure~\ref{perfEva}.
It can be seen that the system can serve more customers, in the presence of a larger percentage of fast charging customers. This is
expected, since the overall ``service rate" is faster in that case.
\begin{figure}[t]
  \includegraphics[width=\columnwidth]{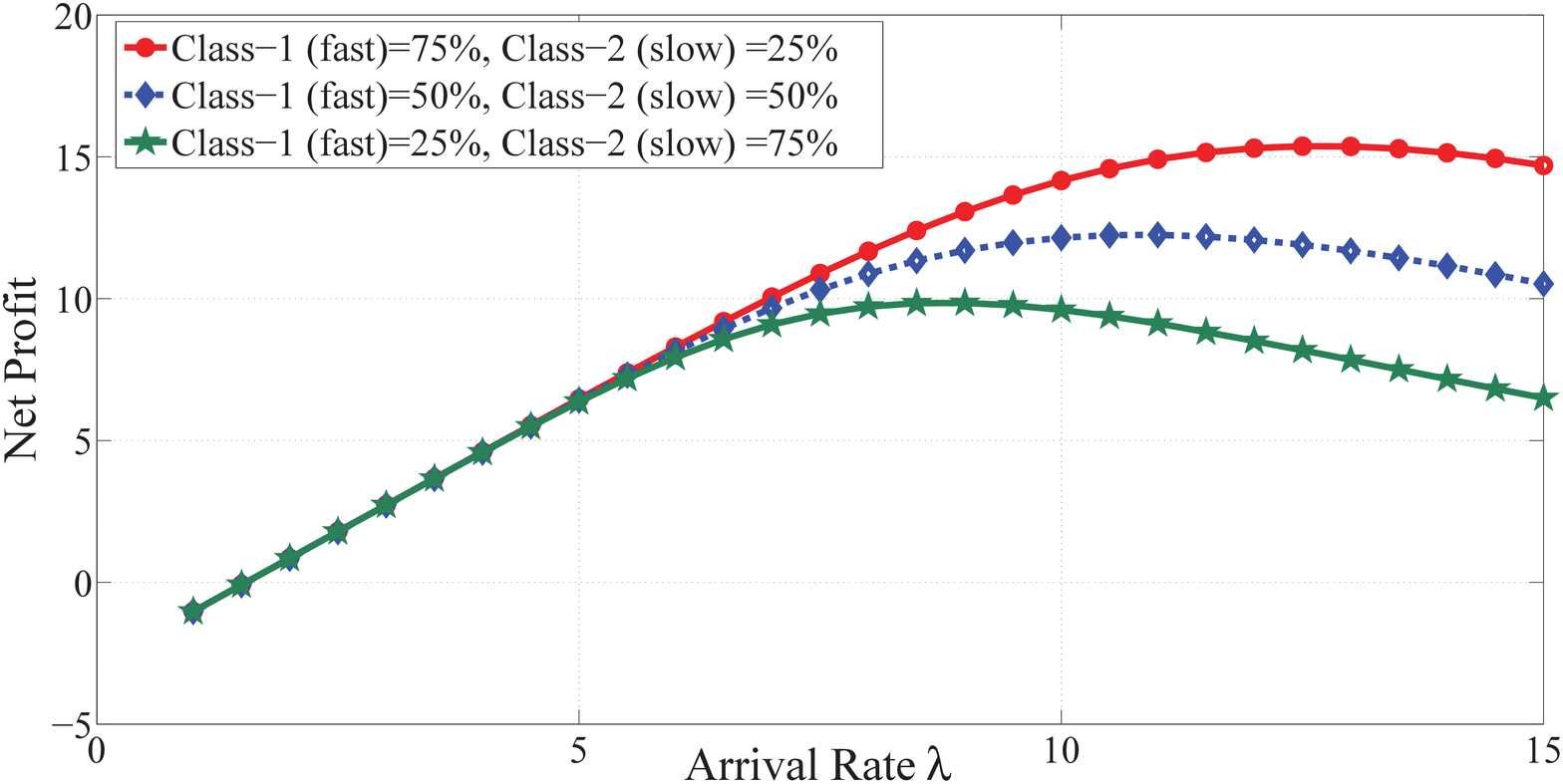}
  \caption{Multi Class Charging Station Net Profit}\label{netProfit}
\end{figure}

\subsection{Profit Model}\label{econModel}
The previous performance assessment provides insight into the gains captured by the posited QoS, namely the blocking probability.
Next, we present a charging station profit model that relates the stochastic model to cost parameters.
This model provides guidance to choose the right values for the amount of power drawn from the grid
for different arrival rates.

The principles of the profit model for $C$ different customer classes are as follows: the charging station earns differential
revenue for each served EV according to its class (e.g. more revenue from fast charging customers etc.). On the other hand, a penalty is paid for each blocked EV because (1) it leads to
dissatisfied customers and degrades the reputation of the station; (2) it enables to control the QoS to foster EV adoption~\cite{ieaRep}; (3) it allows station operators to size its capacity to maximize its profit. It is assumed that a higher penalty is paid to customers charged more
for service. Let $R^{(c)}_g$ and $R^{(c)}_l$ be the revenue gained per EV class-$c$, when served from the grid and the energy storage, respectively. Further, let $C^{(c)}_b$ denote the blocking cost of a single EV in class-$c$. Finally, let $C_0$ represent the fixed installation cost and $C^{(c)}_a R$ the acquisition cost, assumed to be proportional to size, for customer class-$c$ of the storage unit. In order to calculate the net profit, for each customer class, we classify the charging states in the Markov chain
model to: the ``grid charging states" and the ``storage unit charging states". Let $\rho^{\emph{(g)}}=\{(i,j): 0\leq i\leq S, 0\leq j\leq R\}$ denote `the grid charging states" and $\rho^{\emph{(l)}}=\{(i,j): S+j\leq i\leq S+R, 1\leq j\leq R\}$
``the storage unit charging states". Similarly, $\rho^{\emph{(bl)}}$ represents the ``blocking states", while $i(s)$ denotes the number of EVs at state $s$.
Then, the proposed profit function can be written as
\begin{eqnarray}
P=& \sum\limits_{c\in \mathcal{C}} \sum\limits_{s\in\rho^\emph{(g)}} R^{(c)}_g i^{(c)}{(s)}\pi^{(c)}(s) + \sum\limits_{c\in \mathcal{C}}\sum\limits_{s\in\rho^\emph{(l)}} R^{(c)}_l i^{(c)}(s)\pi^{(c)} (s)  \nonumber \\
  & - (C_0+\sum\limits_{c\in \mathcal{C}}R^{(c)} C^{(c)}_a) - \sum\limits_{c\in \mathcal{C}} \sum\limits_{s\in\rho^\emph{(bl)}} C^{(c)}_b i^{(c)}{(s)}\pi^{(c)}(s)\;\;
\end{eqnarray}

We evaluate the profit model for the following set of parameters in the presence of two customer classes (fast/slow):
$R^{(1)}_g$=$R^{(1)}_l$=$3$, $R^{(2)}_g$=$R^{(2)}_l$=$1.5$, $C^{(1)}_b$=$3.5$, $C^{(2)}_b$=$2$,  $C^{(1)}_a$=$0.25$, $C^{(2)}_a$=$0.15$ and $C_0$=$0.02$. The results are shown in Figure~\ref{netProfit}. For low arrival rates, the cost related to acquisition and installation outweighs the revenue gained from charging EVs, and hence a negative profit is earned. On the other hand, for high arrivals rates, the cost of blocking customers becomes dominant and the total net profit decreases. Moreover, since fast charging lowers the blocking probability, the
system means more profit when the proportion of class-1 customers is higher.

\section{A Network of Charging Stations}\label{chargingNetwork}
\subsection{Overview}
Fast public charging stations are key to build confidence in the early stages of EV adoption. At present, the number of fast charging stations in the US is quite low, and deployment plans in the short term are limited to selected highways only~\cite{gigaom,insider}. In order to compete against gas stations, deploying urban charging facilities becomes necessary~\cite{canada1}. In this section, the operation of a network
of fast charging stations in an urban environment is studied, where each individual station is modeled according to the architecture introduced in section~\ref{singleStation}.

In the real world, urban traffic movements are far from being uniform. In fact, people drive between specific points of interest, such as their home, school, workplace, etc. Driving patterns vary according to the time of the day (weekday rush hours, weekends etc.) and hence
traffic density represents a dominant factor in the utilization of each node in a charging station network. As the power grid limitations prevent stations from providing more capacity, grid operators have to consider the fact of spatial and temporal demand to optimally allocate their power resources.

\subsection{Power Resource Allocation in a Charging Station Network }\label{optProblems}
\subsubsection{Case-\rom{1}: No Allocation}\label{noAlloc}
In the first case, all charging stations in the network are assumed to be identical. Let $l = {1,2,...,N}$ be the index set of charging stations. Further assume that each station serves $c\in \mathcal{C}$ types of customer classes, so that $S^{(c)}_1$=$S^{(c)}_2\dots S^{(c)}_N$ and $R^{(c)}_1$ = $R^{(c)}_2 \dots = R^{(c)}_N$. The only parameter that differs in these stations is the arrival rate $\lambda_i$ and
composition of the customer class populations $\vec{\rho}$, which comes from the traffic density (note that we consider rational customers who always drive to the nearest station).
\subsubsection{Case-\rom{2}: Optimal grid power-$S$ allocation within a large geographical urban areas}
Similarly to the case above, there are $N$ charging stations deployed in a large urban environment. However, customers also have access to charging station location information provided by a central authority\footnote{Via smart apps such as~\cite{plugShare} or on board communication systems~\cite{ieeeP2030}}. This case is divided into two subcases. The first subcase assumes that all drivers are selfish and  similarly to Case~\rom{1}, they choose the nearest charging station. The second subcase assumes a hybrid population of selfish drivers and EV fleets. Note that unlike selfish users, EV fleets adhere to the decisions of the power utility to fulfill the requirements of customer agreements. Hence, the arrival rate of each station can be shaped within a $[\lambda_{min},\;\lambda_{max}]$ range.

Let $S_{max}$ be the maximum level of generation capacity that the grid can supply to the network in a metropolitan area. Also each station serves $c\in \mathcal{C}$ types of customer classes. Using the discretization assumption at each charging node, two resource allocation problems are formulated as a mixed integer non-linear programming problem in Equations~\ref{AllocateS} and~\ref{AllocateSwL}. For both subcases, the
proposed scheme allocates more power resources to the busier stations, while taking into account \emph{QoS} targets. If the total power required to satisfy the \emph{QoS} requirements is greater than $S_{max}$, then the charging station network provides best-effort service with the maximum allowable grid power, $S_{max}$.
\paragraph{}
\vspace{-1.0 mm}
\begin{equation}\label{AllocateS}
\begin{aligned}
& \underset{S}\min && \sum\limits_{i\in l} \sum\limits_{c\in \mathcal{C}} B_i(\vec{\rho_i}\lambda_i, S^{(c)}_i, R^{(c)}_i)\\
& \text{s.t. }&&\sum\limits_{i\in l} \sum\limits_{c\in \mathcal{C}} S^{(c)}_i = S\\
&&&0\leq B_i(\vec{\rho}\lambda_i, S^{(c)}_i, R^{(c)}_i) \leq \epsilon  \\
&&& S^{(c)}_i \; \in\mathbb{Z^+}\\
&&& R^{(c)}_i ,\; \lambda_i ,\; and \; \vec{\rho_i},\; are\; given\\
&&&  \forall i\in l,\; \forall c\in \mathcal{C}
\end{aligned}
\end{equation}
\paragraph{}
\begin{equation}\label{AllocateSwL}
\begin{aligned}
& \underset{S,\lambda}\min &&\sum\limits_{i\in l} \sum\limits_{c\in \mathcal{C}} B_i(\vec{\rho_i}\lambda_i, S^{(c)}_i, R^{(c)}_i)\\
& \text{s.t. }&&\sum\limits_{i\in l} \sum\limits_{c\in \mathcal{C}} S^{(c)}_i = S\\
&&& 0\leq B_i(\vec{\rho}\lambda_i, S^{(c)}_i, R^{(c)}_i) \leq \epsilon  \\
&&&\lambda^{(c)}_{min}\leq \lambda^{(c)}_i \leq \lambda^{(c)}_{max}\\
&&& S^{(c)}_i \; \in\mathbb{Z^+}\\
&&& R_i^{(c)}\;, \vec{\rho_i}\;, \lambda^{(c)}_{min}\; and\; \lambda^{(c)}_{max}\; are\; given\\
&&&\forall i\in l,\;\forall c\in \mathcal{C}
\end{aligned}
\end{equation}
\subsubsection{Case-\rom{3}: Optimal $S$ and $\lambda$ allocation in small geographical areas}
Let $l^*\subset l$ and $0<n\leq N$. In this case, a charging station network deployed over a relatively well confined small
geographical area with $n$ stations is considered. This case is different from the previous scenario in the following aspect: the total population consists of EV fleets and through agreements, customers can be assigned to any neighboring station. Since the considered distances between stations are reasonably short ($2-3$ miles\footnote{It would require $0.5-1kWh$ of stored energy and would cost $10$-$20$ cents with the current rates.}), routing customers to other stations would have negligible cost to drivers. Thus, customers can be assigned to neighboring area stations to minimize the total blocking probability. In all cases, the local energy storage is assumed to have already been
acquired by the charging station(e.g. $R=5$), thus its size is fixed. Finally, each station serves $c\in \mathcal{C}$ classes of customers, and routed customers get the same type of service. Then, the optimization problem becomes:
\begin{equation}\label{SLAllocation}
\begin{aligned}
& \underset{S,\lambda}\min &&\sum\limits_{i\in l} \sum\limits_{c\in \mathcal{C}} B_i(\vec{\rho_i}\lambda_i, S^{(c)}_i, R^{(c)}_i)\\
& \text{s.t. }  &&\sum\limits_{i\in l} \sum\limits_{c\in \mathcal{C}} S^{(c)}_i = S\\
&&& \sum\limits_{i\in l^*} \lambda_i = \lambda\\
&&& 0\leq B_i(\vec{\rho_i}\lambda_i, S^{(c)}_i, R^{(c)}_i) \leq \epsilon  \\
&&& \lambda^{(c)}_i \geq 0,\; \forall i\in l^*\\
&&& S^{(c)}_i\in \mathbb{Z^+}\\
&&& R_i\;and\;\vec{\rho_i}\; are\; given\\
&&& \forall i\in l^*,\;\forall c\in \mathcal{C}
\end{aligned}
\end{equation}
In addition to the system constraints presented for each allocation problem, there may be additional constraints, depending on the existing power network, such as distribution network limitations, etc. However, since the interaction of the charging stations with the grid is limited to the constant power drawn from it, these case-by-case varying constraints may only affect the maximum power allocation for individual stations. Thus, these constraints can easily be incorporated within the existing formulation to address these allocation problems.
\begin{figure}[t]
  \centering
\includegraphics[width=0.75\columnwidth]{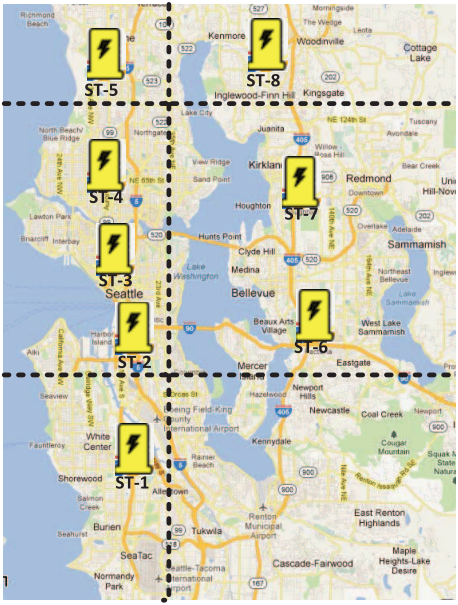}
  \caption{Fast DC charging station map in Seattle, WA ~\cite{adhoc}}\label{seattle}
  \end{figure}

\section{Evaluation \& Results}\label{results}
\subsection{Overview}
Collecting vehicular traffic traces, especially in urban areas, is a challenging and costly task. Hence, vehicles movements are not well
calibrated. However, in~\cite{adhoc} bus movements from the Seattle area were obtained. Due to the city's physical layout and
extensive bus network\footnote{$1200$ buses in a \~5000 square kilometers area}, it is claimed that these movements resemble actual traffic patterns quite closely. In the next subsection (\ref{inputAnalysis}), we use this publicly available data to investigate the spatial distribution of vehicles, during weekday rush hour (7am-9am and 5pm-7pm). The remainder of this subsection is organized as follows. In subsection~\ref{StationRep}, we explain our methodology in locating fast charging stations on the city map. In subsection~\ref{metamodel}, we use the Response Surface Methodology to approximate charging station blocking probabilities into a second order regression metamodel. Finally, in subsection~\ref{comparisons}, we solve the optimization problems presented in section~\ref{optProblems} using our metamodel.
\subsubsection{Input Analysis}\label{inputAnalysis}
According to~\cite{adhoc}, the location of each bus was recorded frequently. We start by normalizing the $x$ and $y$ coordinates of the input data.
Subsequently, the ARENA Input Analyzer~\cite{arena} is used to fit a spatial distribution to the data. The results indicated that with
mean squared error of $0.6\%$, the spatial distribution of vehicles is a piecewise beta distribution for weekday rush hours. The results are presented in equations~\ref{spatDistX} and~\ref{spatDistY}.
\begin{equation}\label{spatDistX}
f(X) =
\begin{cases}
    44\times BETA(4.42, 0.763) & \text{ 0 $\le$ X $\le$ 44} \\
    44 + 137\times BETA(0.752, 4.7) & \text{44 $\le$ X $\le$ 180}
\end{cases}
\end{equation}
\begin{equation}\label{spatDistY}
f(Y) =
\begin{cases}
    150\times BETA(2.42, 0.799) & \text{ 0 $\le$ Y $\le$ 150} \\
    150 + 121\times BETA(1.07, 5.44) & \text{150 $\le$ Y $\le$ 270}
\end{cases}
\end{equation}

In addition, we analyze the correlation of $x$ and $y$ coordinates, and calculate the correlation coefficient as $0.06$.

\begin{figure}[t]
  \centering
\includegraphics[width=\columnwidth]{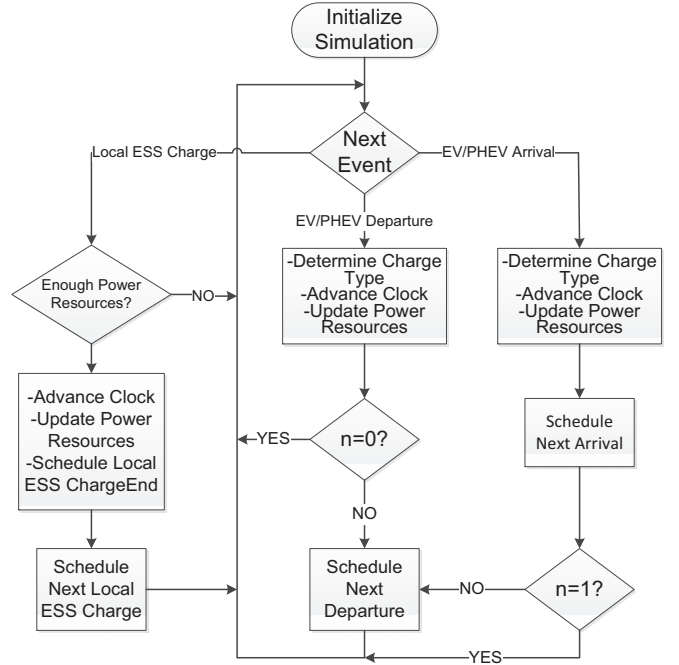}
  \caption{Discrete Event Simulation Flow Chart}\label{fChart}
  \end{figure}
\subsubsection{Charging Station Placement}\label{StationRep}
In~\cite{adhoc}, researchers placed eight base station towers in such a way that base stations can communicate with all mobile nodes. Since the charging station layout problem is outside the scope of this paper, a similar approach is used and the same number of charging stations is
deployed in the same locations given by the following coordinates: $\{x_i,y_i\}=\{60,45\}\bc \{60,90\}\bc \{60,135\}\bc \{60,180\}\bc \{60,225\}\bc \{100,90\}\bc \{100,160\}\bc \{100,225\}$. Figure~\ref{seattle} presents the map with the locations of the charging stations. In order to calculate the traffic intensity at each station, a discrete event simulation model is used. We present its flowchart in Figure~\ref{fChart}. The station parameters are given by $S=5$, $R=5$, $\mu=2$, $\nu=4$ (assuming only fast charging customers) for all stations. The simulation is terminated when one million vehicles get serviced. It is run for a total of $30$ times and $95\%$ confidence intervals of the parameters of interest are
obtained. The traffic intensity for each station is shown in Table~\ref{arrivalRates}.
It can be seen from Table~\ref{arrivalRates} that charging stations three and four are used to meet most of the charging demand, whereas other stations have relatively little demand. For instance, letting the overall arrival rate be $\lambda=50$, then the blocking probabilities for the eight identical  stations would be \begin{small}$\vec{B_i}=[0.019,\;0.053,\;0.58,\;0.58,\;0.0158,\;0.0153,\;0.043,0.014]$\end{small}. It can be concluded from this expository calculation that, without any power allocation, there could be severe fluctuations in terms of \emph{QoS} among the charging facilities\footnote{Some stations(e.g. station-$3$) will exhibit a very high blocking probability, whereas overprovisioned stations (e.g. station-$8$) will exhibit a very low blocking probability.}.

\subsection{Output Analysis}\label{output}
\subsubsection{Metamodeling of Blocking Probabilities}\label{metamodel}
\begin{table}
  \centering
  \caption{Traffic Intensity (T.I.) of Each Station} \label{arrivalRates}
    \begin{tabular}{c|c|c|c|c|c}
    \hline
    Sta. ID & mean(T.I.) & 95\% CI & Sta. ID & mean(T.I.) & 95\% CI \\ \hline\hline
    1     & 3.56\% & 0.030\% & 2     & 9.68\%& 0.022\% \\
    3     & 36.99\% & 0.060\% & 4     & 36.68\% & 0.035\% \\
    5     & 1.85\% & 0.037\% & 6     & 1.65\% & 0.018\% \\
    7     & 8.6\% & 0.055\%   &  8     & 0.98\%& 0.016\% \\
    \hline
    \end{tabular}%
\end{table}%
In section~\ref{singleStation}, numerical methods are used to calculate EV blocking probabilities. However, new calculations are needed for
each set of new input parameters $S$, $R$, $\nu$, and $\lambda$ to determine the blocking probability $B$. using the
 Response Surface Methodology (RSM) we are able to calculate an approximate second order polynomial model for the functional relationship
 between $B$ and the input parameters ($B = f(S,R,\nu,\lambda)$~\cite{rsmBook}. As input parameters we used
 those presented in Table~\ref{rsmTable}, keeping $\mu$ is fixed to $2$.
\begin{table}[h]
\centering
\caption{RSM Input Parameters}\label{rsmTable}
\begin{tabular}{c|c|c|c}
  \hline
  Parameter & Interval & Increments & Type \\
  \hline\hline
  S & [1,15] & 1 & Integer \\ \hline
  R & [1,15] & 1 & Integer \\ \hline
  $\lambda$ & [0.25,30] & 0.25 & Float \\ \hline
  $\nu$ & [2,10] & 1 & Integer \\ \hline
\end{tabular}
\end{table}
The handicap of this approach is that the blocking probabilities have to be in the $[0,1]$ interval, whereas the RSM model
can predict values outside it. For that reason, we fit the RSM model to the logit transformation ($ y = log(x/(1 - x))$)
of $B$ and then use the inverse-logit ($x = 1/(1+e^{-y})$) transformation to obtain the final results.
The regression model where the response variable corresponds to logit($B$) is given in Equation~\ref{rsmEqu}.

\begin{align}\label{rsmEqu}
 B(S, R, \lambda, \nu)=-3.990-2.666S-1.6152R-0.1492\nu\nonumber \\
   +3.840\lambda-0.0645SR-0.002S\nu+0.209S\lambda-0.0078R\nu\nonumber  \\
   +0.094R\lambda + 0.003\nu\lambda -0.0175S^2+0.055R^2+0.0089\nu^2\nonumber \\
     -0.271\lambda^2
\end{align}
Then, the blocking probability becomes,
\begin{equation}
 Blocking\; Prob.= \left\{ {\begin{array}{*{20}{c}}
{B( \cdot )}\\
0
\end{array}} \right.\begin{array}{*{20}{c}}
{if\lambda  > 0}\\
{if\lambda  = 0}
\end{array}
\end{equation}

For the above regression model, the R-Square statistic is $88.06\%$ and the mean square root error is $0.52\%$.
\begin{figure}[h]
\begin{eqnarray}\label{jacobian}
\setlength{\extrarowheight}{3.5pt}
\begin{bmatrix}
 \partial B\over\partial S  \\
 \partial B\over\partial R\\
  \partial B\over\partial \nu\\
 \partial B\over\partial \lambda
\end{bmatrix}
=
\begin{small}
\begin{bmatrix}

    -0.035S-0.0645R-0.002\nu+0.21\lambda-2.66 \\
    0.0014S+0.11R-0.008\nu+0.094\lambda-1.62 \\
    -0.002S-0.078R+0.178\nu+0.025\lambda-0.15 \\
    0.209S+0.094R+0.003\nu-0.54\lambda+3.84

\end{bmatrix}
    \end{small}
\end{eqnarray}
\end{figure}

\begin{figure}[h]
\begin{eqnarray}\label{Hessian}
\emph{H}
=
\begin{bmatrix}
 -0.035 &-0.0645 &-0.002 &0.21 \\
    0.0014 &0.11 &-0.008 & 0.094 \\
    -0.002 &-0.078 &0.178& 0.025 \\
    0.209 &0.094 &0.003 &-0.54
\end{bmatrix}
\end{eqnarray}
\end{figure}

Some key quantities like the Jacobian (equation \ref{jacobian}) and the Hessian matrix (equaiton \ref{Hessian}) are given to aid assessing the sensitivity of $B$
with respect to inputs $S$, $R$, $\nu$ and $\lambda$ variables is presented. It can be seen that grid power $S$ has the highest impact for decreasing the blocking probability. Note that during periods of high arrival rates, there is going to be little spare capacity left and
hence the local storage device would be frequently in an empty state, as indicated by these results.
\begin{table*}[t]
\centering
\caption{Results for Case~\rom{2}B (Mixed Population of Selfish EVs and Fleets)}\label{case2-2}
    \begin{tabular}{cc |ccc|ccc|ccc|ccc}
    \hline
          &   $\sum\limits_{i\in l} \lambda_i $    & \multicolumn{3}{c|}{$Station_1$} & \multicolumn{3}{c|}{$Station_2$} & \multicolumn{3}{c|}{$Station_3$} & \multicolumn{3}{c}{$Station_4$} \bigstrut\\
    \hline\hline
          &  & $S_1$    & $\lambda_1$    & $B_1$    & $S_2$    & $\lambda_2$    & $B_2$    & $S_3$    & $\lambda_3$    & $B_3$    & $S_4$    & $\lambda_4$    & $B_4$ \bigstrut \\
    \multicolumn{1}{c}{\multirow{3}[0]{*}{\begin{sideways} {$\epsilon$=0.05} \end{sideways}}} & 20    & 1     & 0.77 & 0.0037 & 2     & 2.656 & 0.015 & 5     & 7.3 & 0.0347 & 5     & 6.984 & 0.032 \\
    \multicolumn{1}{c}{} & 25    & 1     & 0.9625 & 0.0153 & 3     & 2.832 & 0.0125 & 6    & 8.98 & 0.047 & 6    & 8.90 & 0.041 \\
    \multicolumn{1}{c}{} & 30    & 1     & 1.155 & 0.032 & 3     & 3.3984 & 0.031 & 7    & 10.767 & 0.05 & 7    & 10.4760 & 0.046 \\ \hline
    \multicolumn{1}{c}{\multirow{3}[0]{*}{\begin{sideways} {$\epsilon$=0.10} \end{sideways}}} & 20   & 1     & 0.77 & 0.057 & 1     & 2.2656 & 0.0317 & 4     & 7.178 & 0.10  & 4     & 6.98 & 0.0932 \\
    \multicolumn{1}{c}{} & 25    & 1     & 0.9625 & 0.0153 & 1     & 2.832 & 0.03491 & 5     & 8.9725 & 0.10 & 5     & 8.73 & 0.0936 \\
    \multicolumn{1}{c}{} & 30    & 1     & 1.155 & 0.032 & 2     & 3.3984 & 0.031 & 6    & 10.76 & 0.10 & 6    & 10.476 & 0.0934 \\
          \hline\hline
          & $\sum\limits_{i\in l} \lambda_i $     & \multicolumn{3}{c|}{$Station_5$} & \multicolumn{3}{c|}{$Station_6$} & \multicolumn{3}{c|}{$Station_7$} & \multicolumn{3}{c}{$Station_8$} \bigstrut \\
          \hline
          &  & $S_5$    & $\lambda_5$    & $B_5$    & $S_6$    & $\lambda_6$    & $B_6$    & $S_7$    & $\lambda_7$    & $B_7$    & $S_8$    & $\lambda_8$    & $B_8$ \bigstrut \\
    \multicolumn{1}{c}{\multirow{3}[0]{*}{\begin{sideways} {$\epsilon$=0.05} \end{sideways}}} & 20    & 1     & 0.3960 & 0.00122 & 1     & 0.3663 & 0.0007 & 2     & 1.892 & 0.0087 & 1     & 0.022 & 0 \\
    \multicolumn{1}{c}{} & 25    & 1     & 0.495 & 0.0007 & 1     & 0.4538 & 0.0004 & 2     & 2.365 & 0.0379 & 1     & 0.0275  & 0 \\
    \multicolumn{1}{c}{} & 30    & 1     & 0.5940 & 0.0017 & 1     & 0.7445 & 0.0011 & 3     & 2.84 & 0.0127 & 1     & 0.033 & 0 \\ \hline
    \multicolumn{1}{c}{\multirow{3}[0]{*}{\begin{sideways} {$\epsilon$=0.10} \end{sideways}}} & 20    & 1     & 0.396 & 0.0002 & 1     & 0.3630 & 0.0001 & 2     & 1.892 & 0.014 & 1     & 0.022 & 0 \\
    \multicolumn{1}{c}{} & 25    & 1     & 0.795 & 0.0007 & 1     & 0.4538 & 0.0004 & 2     & 2.365 & 0.0379 & 1     & 0.1875 & 0 \\
    \multicolumn{1}{c}{} & 30    & 1     & 0.594 & 0.0017 & 1     & 05445 & 0.0011 & 2     & 2.883  & 0.0764 & 1     & 0.651 & 0.001 \\
    \hline
    \end{tabular}%
  \label{tab:addlabel}%
\end{table*}%
\subsection{Comparison of three cases}\label{comparisons}

\begin{figure}[t]
  \centering
\includegraphics[width=\columnwidth]{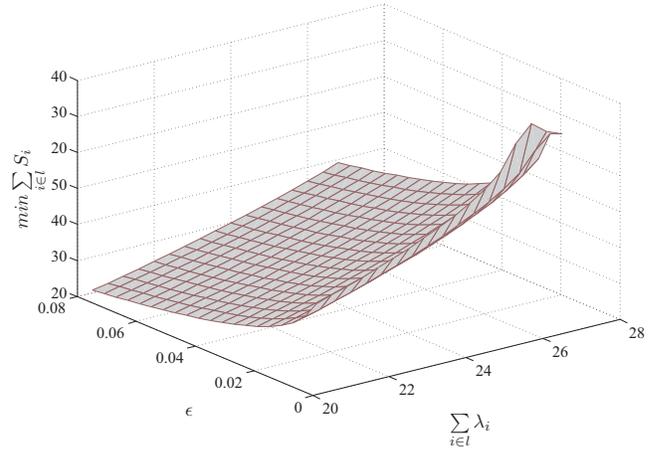}
  \caption{Minimum grid power required to meet $\epsilon$ \emph{QoS} targets}\label{minSS}
  \end{figure}
\setlength{\tabcolsep}{4pt}
\begin{table}
  \centering
  \caption{Results for Case-\rom{2}A (Selfish EVs)}\label{case2-1}
    \begin{tabular}{ll|cc|cc|cc|cc}%
    \hline
          &   $\sum\limits_{i\in l} \lambda_i $    & \multicolumn{2}{c|}{$Station_1$} & \multicolumn{2}{c|}{$Station_2$} & \multicolumn{2}{c|}{$Station_3$} & \multicolumn{2}{c}{$Station_4$}  \bigstrut \\
    \hline\hline
    \multicolumn{1}{c}{\multirow{5}[0]{*}{\begin{sideways} {$\epsilon$=0.05} \end{sideways}}} &       & $S_1$    & $B_1$    & $S_2$    & $B_2$    & $S_3$    & $B_3$    & $S_4$    & $B_4$ \bigstrut \\
    \multicolumn{1}{c}{} & 20  & 1     & 0.0023 & 2     & 0.0094 & 6    & 0.025 & 6    & 0.0214 \\
    \multicolumn{1}{c}{} & 25  & 1     & 0.0067 & 2     & 0.0269 & 7    & 0.0348  & 7    & 0.0299 \\
    \multicolumn{1}{c}{} & 30  & 1     & 0.0148 & 3     & 0.0082 & 8    & 0.0436 & 8   & 0.0376 \\ \hline
    \multicolumn{1}{c}{\multirow{3}[0]{*}{\begin{sideways} {$\epsilon$=0.10} \end{sideways}}} & 20  & 1     & 0.023 & 2     & 0.0094 & 5  & 0.0641 & 5     & 0.0568  \\
    \multicolumn{1}{c}{} & 25  & 1     & 0.0067 & 2     & 0.0269 & 6    & 0.0754 & 6    & 0.0669 \\
    \multicolumn{1}{c}{} & 30  & 1     & 0.0148 & 2     & 0.0568 & 7    & 0.0842 & 7    & 0.0748 \\
    \hline
          &   $\sum\limits_{i\in l} \lambda_i $   & \multicolumn{2}{c|}{$Station_5$} & \multicolumn{2}{c|}{$Station_6$} & \multicolumn{2}{c|}{$Station_7$} & \multicolumn{2}{c}{$Station_8$} \bigstrut \\
    \hline\hline
    \multicolumn{1}{c}{\multirow{4}[0]{*}{\begin{sideways} {$\epsilon$=0.05} \end{sideways}}} &       & $S_5$    & $B_5$    & $S_6$    & $B_6$    & $S_7$    & $B_7$    & $S_8$    & $B_8$ \bigstrut \\
    \multicolumn{1}{c}{} & 20  & 1     & 0.0001 & 1     & 0 & 2     & 0.0054 & 1     & 0 \\
    \multicolumn{1}{c}{} & 25  & 1     & 0.0002 & 1     & 0.0002 & 2     & 0.063 & 1     & 0 \\
    \multicolumn{1}{c}{} & 30  & 1     & 0.0006 & 1     & 0.0004 & 2     & 0.0367 & 1     & 0 \\ \hline
    \multicolumn{1}{c}{\multirow{3}[0]{*}{\begin{sideways} {$\epsilon$=0.10} \end{sideways}}}  & 20  & 1     & 0.0001 & 1     & 0 & 2     & 0.0911 & 1     & 0  \\
    \multicolumn{1}{c}{} & 25  & 1     & 0.0002 & 1     & 0.0002 & 2     & 0.0163 & 1     & 0 \\
    \multicolumn{1}{c}{} & 30  & 1     & 0.0006 & 1     & 0.0004 & 2     & 0.0367 & 1     & 0 \\
    \hline
    \end{tabular}%
\end{table}%
  \begin{table}
    \centering  \caption{Comparison of Case-\rom{2}A (Selfish EVs) and Case-\rom{2}B (Mixed Population)}\label{case2-3}\label{case2comp}
    \begin{tabular}{cc|c|c|c|r}
    \hline
          & & $\sum\limits_{i\in l} \lambda_i$  & $S^{Case~\rom{2}A}_i $ & $ S^{Case~\rom{2}B}_i $ & $Savings$   \bigstrut \\
    \hline\hline
    \multicolumn{1}{c}{}{\multirow{3}[0]{*}{\begin{sideways}{$\epsilon$=0.05}\end{sideways}}} & & 20 & 20& 18 & 10\%\\
    \multicolumn{1}{c}{} & & 25  & 22 & 21 &  4.55\%\\
    \multicolumn{1}{c}{} & &30  & 25 & 24&  4\%    \\
    \hline
    \multicolumn{1}{c}{\multirow{3}[0]{*}{\begin{sideways}{$\epsilon$=0.10}\end{sideways}}}& & 20 & 18&15& 18.75\%\\
    \multicolumn{1}{c}{} & & 25  & 20 & 17 & 15\% \\
    \multicolumn{1}{c}{}& & 30  & 22 & 20 &9\%   \\
    \hline
\end{tabular}%
\end{table}%
\setlength{\tabcolsep}{2pt}
\begin{table}
  \centering
  \caption{Results for Case~\rom{3} (EV Fleets)}\label{case3}
    \begin{tabular}{cc|ccc|ccc|ccc}
    \hline
      $\sum\limits_{i\in l^*} S_i$   & $\sum\limits_{i\in l^*} \lambda_i$ & \multicolumn{3}{c|}{$Station_2$} & \multicolumn{3}{c|}{$Station_3$} & \multicolumn{3}{c}{$Station_4$} \\
    \hline\hline
    & & $S_2$    & $\lambda_2$    & $B_2$    & $S_3$    & $\lambda_3$    & $B_3$    & $S_4$    & $\lambda_4$    & $B_4$ \bigstrut\\
    18 & 16.67 & 6     & 5.56 & 0.004& 6     & 5.56 & 0.004 & 6     & 5.56 & 0.004 \\
    24 & 20.83 & 8     & 6.94 & 0.0327 & 8     & 6.94 & 0.0327 & 8     & 6.94 & 0.0327 \\
    30 & 24.9 & 10    & 8.3 & 0.005 & 10    & 8.3 & 0.005 & 10    & 8.3 & 0.005 \\
    \hline
    \end{tabular}%
\end{table}%
Next, we compare the performance of the following three scenarios:
\begin{inparaenum}[(i)]
\item all eight stations are identical (case-\rom{1});
\item power resource allocation for selfish EV population (case-\rom{2}A) and mixed (selfish and EV fleets) population (case-\rom{2}B); and
\item power resource allocation for EV fleets only (case-\rom{3}).
\end{inparaenum}
 Standard interior point methods are used to solve the optimization problems introduced in section~\ref{optProblems}. Problems formulated in case-\rom{2} and case-\rom{3} are non-linear integer programs and they are solved by relaxing the integer constraint and ceiling to the nearest integer. For case-\rom{2}A, suppose that the station operator wants to provide $\epsilon$-level \emph{QoS} at all stations. One of the main goals of this scheme is to use the minimum amount of power grid resources (for illustration assume all customers demand fast charging). Hence, the minimum required grid power $S_{min}$ to meet the~\emph{QoS} targets is calculated. As long as $S_{min} \leq S_{max}$ where $S_{max}$ is the total allocated generation capacity, this target is going to be reached. A generic calculation is presented in Figure~\ref{minSS}.

Next, let us compare cases-\rom{1} and -\rom{2}A. Suppose that the charging station operator wants to ensure that each station can meet $90\%$ of the customer demand at all times ($\epsilon$=$0.10$). For the eight stations, the arrival rate is assumed to be $\lambda$ = $27$. Since the majority of the population resides near Stations $2$, $3$ and $4$, we assume that these two stations serve two types of customers; class-$1$
(fast charging $\mu$=$2$) and class-$2$ (slow charging $\mu$=$1$). The same set of parameters from section~\ref{queueingModel} are used for the efficiency and the power rating of the local energy storage units. Since these regions are close to downtown we further assume that $\vec{\rho}$=$(75\%,25\%)$. The remainder of the stations serve customer class-$1$. Solving equation~\ref{AllocateS} results in $\vec{S}$=[1, 2, 9, 9, 1, 1, 2, 1]. With the allocated grid power, blocking probabilities for each station are $\vec{B}$=$[0.0094, 0.028, 0.099, 0.087, 0.0004, 0.0023, 0.016, 0.0001]$. In order to compare the performance of the whole network, we calculate the weighted sum of stations' blocking probabilities;
\begin{equation}
\sum\limits_{i \in l} {{w_i}{B_i}},\; \text{where}\;{w_i} = \frac{{{\lambda _i}}}{{\sum\limits_{i \in l} {{\lambda _i}} }}
\end{equation}
Then, the weighted sum of blocking becomes $\sum\limits_{i \in l} {{w_i}{B_i}}$=$0.0440$. To compare these results with case \rom{1}, assume that each station has $S_i=3$ (except $S_3=S_4=4$) and $R=5$. Arrival rates are the same as case \rom{2}. For this case, the weighted sum of blocking becomes $\sum\limits_{i \in l} {{w_i}{B_i}}$=$0.4365$. This sample calculation shows that with power resource allocation, more vehicles can receive service with the same amount of grid power.
\begin{figure}[t]
\centering
\includegraphics[width=\columnwidth]{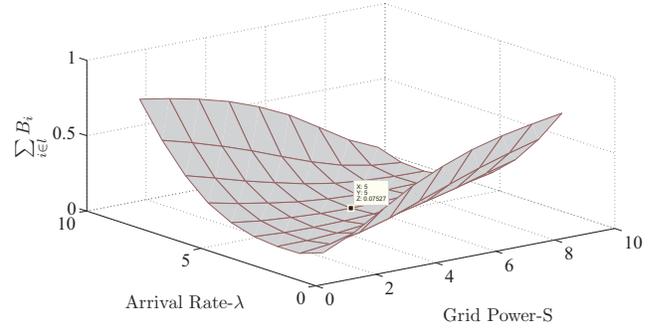}
\caption{Evaluation of Equation~\ref{SLAllocation}}\label{SLAllocationResults}
\end{figure}

For the allocation problem in cases~\rom{2}-B and ~\rom{3}, a two-way communication infrastructure is used to offer customers incentives to charge from other stations. In the first case, a central authority can route a certain percentage of customers in the $[\lambda_{min}, \lambda_{max}]$ range. In the latter one, any customer can be assigned to any station in the same neighboring area. Thus, for the first case suppose that the arrival rate of each station is in the $\lambda\pm10\%$ range. For instance, if station-$3$ the arrival rate is $\lambda_3=10$ arrival rate will be assigned in [9,11] interval. Next, two allocation problems~\ref{AllocateS} and~\ref{AllocateSwL} are solved for six different combinations: $\epsilon=0.05,\;0.10$ and $\sum\limits_{i\in l} \lambda_i = 20,\;25,$ and $30$.

We proceed to compare three cases for all stations in detail. Table~\ref{case2-1} presents the results for a population of selfish EV users. The utility can only allocate optimal power (problem~\ref{AllocateS}). On the other hand, for a mixed population of EVs (selfish and Fleets) allocation problem~\ref{AllocateSwL} is solved. We present detailed results for each charging facility in Table~\ref{case2-2}. Note that since the central authority can partially affect the customer choices, blocking probability targets can be achieved with less grid power. For instance, customer routing can lead to $10\%$ power savings to provide $\epsilon$=$0.05$ QoS. Detailed results are given in Table~\ref{case2-3}. Moreover, the comparison of these two cases and associated savings are presented in Table~\ref{case2comp}.
\begin{figure}[t]
  \centering
\includegraphics[width=\columnwidth]{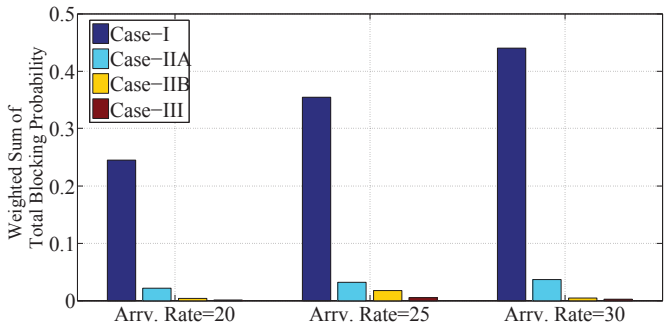}
  \caption{Comparison of three cases respect to stations 2, 3 and 4}\label{barAll}
  \end{figure}

Case-\rom{3} assumes a population of pure EV fleets. Note that the network map (Figure~\ref{seattle}) is divided into smaller geographic areas, and inside each region the cost of driving between charging stations has a negligible cost. We assume that charging stations $2$, $3$, and $4$ constitute a charging network. Similar to the previous case, a central authority through the use of two-way communications, can assign customers to any station in this subarea. To evaluate this case, assume that stations $2$, $3$, and $4$ are in a small well-confined neighborhood and driving between these stations has a negligible cost. Then, in optimization problem~\ref{SLAllocation}, minimum customer blocking probabilities are obtained at $S_i$=$S/N$ and $\lambda^{(1)}_i$=$\lambda^{(1)}/N$ where $N$ is the number of charging stations (also $R^{(1)}_1=R^{(1)}_2\dots R^{(1)}_N$ and $R^{(2)}_1=R^{(2)}_2\dots R^{(2)}_N$)\footnote{Note that customer profile $\vec{\rho}$ is approximately the same since stations are physically close to each other}. We present the results for $\epsilon=0.05$ and varying arrival rate parameters in Table~\ref{case3}.
Moreover,  we run a sample calculation for two stations, with the following parameters: $\sum\limits_{i\in l} \lambda_i = 10$, $\sum\limits_{i\in l} S_i = 10$ and $R_1=R_2=5$. The results are shown in Figure~\ref{SLAllocationResults}.

\begin{figure}[t]
  \centering
\includegraphics[width=\columnwidth]{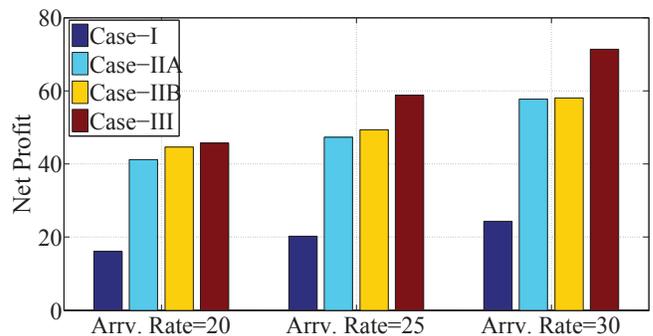}
  \caption{Net Profit Comparison}\label{costComp}
  \end{figure}

Next, we compare the baseline scenario (no allocation of any kind) and three allocation schemes for stations $2$, $3$, and $4$ since they serve both fast and slow customer classes. For a fair comparison, we fixed the total grid resources, and employ the same type of energy storage devices. Results are depicted in Figure~\ref{barAll} for three different arrival rates.
In order to quantify the effects of power allocation and customer routing on the charging network, the profit model of section~\ref{econModel} is applied to all stations. Previously presented results for all three cases are used for $\epsilon = 0.05$ and arrival rates $\lambda=20, \;25,$ and $30$. The same set of parameters from section~\ref{econModel} is employed. In Figure~\ref{costComp}, average net profit per charging station is depicted, which shows that the proposed framework improves both the system (in terms of QoS) and its financial performance significantly.
\section{Towards a more realistic model: the role of communications and incentives}

In this study, we have presented the architecture and an associated stochastic model for a network of charging stations
for allocating power and reroute customers in an optimal fashion. Deployment of such a network requires the necessary
communications infrastructure and protocols to ensure timely information dissemination.

Hence, it is important to quantify the impact of communications on the network's operations. For example,
how communication delays and losses between EVs, charging stations and a network coordinator affect routing decisions and hence QoS?

Wireless network technology such as 4G (LTE or WiMAX) and 3G/UMTS, could play a critical role in EV/PHEV roaming schemes. Moreover, in potential ``handoff" circumstances, the communication network should support inter-grid communications, so that drivers can retrieve up-to-date information about nearest charging stations, available pricing, etc., through the aforementioned smart apps or on-board energy management systems~\cite{ieeeP2030}. In the latter case, charging stations will be in constant communication with the network coordinator, possibly using a mixture of standard wide area communication protocols. In a recent real world application, charging stations employed a 2.4 GHz, 802.15.4 full mesh radio protocol~\cite{wPaper}. Moreover, depending on the volume and the criticality of operations, IEC 60870-6 (inter-control center communications) could be employed as well.

The major problem in information dissemination is related to the network connectivity, that is, some vehicles could be missing part of the information since they are temporarily not connected to the network. Assuming that information dissemination is provided by means of a specific network, this aspect can be translated into the problem of evaluating the degree of connectivity of the employed communication network.

Another limitation of the proposed framework is that it assumes that EVs strictly adhere to routing decisions taken by the network coordinator.
However, many drivers may deviate from the proposed assignments by the network coordinator and head towards the nearest station.
Pricing incentives could address this issues, as discussed in ~\cite{dBan}.

\ifCLASSOPTIONcaptionsoff
  \newpage
\fi


%

%
\bibliographystyle{IEEEtran}
\bibliography{jsac}

\end{document}